\documentclass [prd,eqsecnum,twocolumn,amsfonts,amssymb]{revtex4}

\input epsf

\usepackage[usenames]{color}
\usepackage{graphicx}
\usepackage{bm}
\usepackage{rotating}

\newcommand{\beq}{\begin{equation}}
\newcommand{\eeq}{\end{equation}}
\newcommand{\beqs}{\begin{eqnarray}}
\newcommand{\eeqs}{\end{eqnarray}}

\begin{document}

\title{Exponential Growth Constants for Spanning Forests on Archimedean
  Lattices: Values and Comparisons of Upper Bounds} 

\author{Shu-Chiuan Chang$^a$ and Robert Shrock$^b$}

\affiliation{(a) \ Department of Physics, National Cheng Kung University,
Tainan 70101, Taiwan} 

\affiliation{(b) \ C. N. Yang Institute for Theoretical Physics and
Department of Physics and Astronomy \\
Stony Brook University, Stony Brook, NY 11794, USA }

\begin{abstract}

  We compare our upper bounds on the exponential growth constant
  $\phi(\Lambda)$ characterizing the asymptotic behavior of spanning
  forests on Archimedean lattices $\Lambda$ with recently derived
  upper bounds.  Our upper bounds on $\phi(\Lambda)$, which are very
  close to the respective values of $\phi(\Lambda)$ that we have
  calculated, are shown to be significantly better for these lattices
  than the new upper bounds.

\end{abstract}

\maketitle

\pagestyle{plain}
\pagenumbering{arabic}


\section{Introduction}
\label{intro_section}

Let $G=(V,E)$ be a graph, defined by its vertex and edge sets $V$ and
$E$, and denote $n=n(G)=|V|$ and $e(G)=|E|$ as the numbers of vertices
(= sites) and edges (= bonds) in $G$.  An important problem in
mathematics is the determination of the number of subgraphs of $G$
that satisfy a specified property, and, in particular, the asymptotic
behavior of this number as $n(G) \to \infty$.  A spanning subgraph of
$G$, denoted $G'$, is a graph with the same vertex set, $V'=V$ and a
subset of the edge set of $G$, $E' \subseteq E$.  A forest is a
spanning subgraph that does not contain any cycles. Given a graph $G$,
let us denote the number of spanning forests in $G$ as $N_{SF}(G)$ and
the number of connected spanning subgraphs in $G$ as $N_{CSSG}(G)$. 
For many families of graphs $G$, $N_{SF}(G)$ and $N_{CSSG}(G)$ grow
exponentially rapidly as functions of $n(G)$ for large $n(G)$, 
thereby motivating the definitions of corresponding exponential
growth constants $\phi(\{ G \})$ and $\sigma(\{ G \})$, 
\beq
\phi( \{ G \} ) = \lim_{n(G) \to \infty} [N_{SF}(G)]^{1/n(G)}
\label{phi}
\eeq
and
\beq
\sigma( \{ G \} ) = \lim_{n(G) \to \infty} [N_{CSSG}(G)]^{1/n(G)} \ , 
\label{sigma}
\eeq
where $\{ G \}$ denotes the $n(G) \to \infty$ limit of the graphs in a
given family. Recall that the degree $\Delta_{v_i}$ of a vertex $v_i$
in a graph $G$ is the number of edges connecting to $v_i$.  A graph
with the property that all of its vertices have the same degree
$\Delta$ is termed a $\Delta$-regular graph.  To avoid unimportant
complications, we restrict here to loopless graphs.

In Ref. \cite{sfc} we calculated upper bounds on $\phi(\Lambda)$ and
$\sigma(\Lambda)$, where $\Lambda$ denotes the $n(G) \to \infty$ limit of
an Archimedean lattice graph.  Here an Archimedean lattice is defined
as a uniform tiling of the plane with one or more types of regular
polygons, such that all vertices are equivalent, and hence is
$\Delta$-regular.  In general, an Archimedean lattice $\Lambda$ is
identified by the ordered sequence of regular polygons traversed in a
circuit around any vertex \cite{gsbook,wn}:
\beq
\Lambda = (\prod p_i^{a_i}) \ , 
\label{arch}
\eeq
where the $i$'th polygon has $p_i$ sides and appears $a_i$ times
contiguously in the sequence (it can also occur
non-contiguously). There are three Archimedean lattices which each
involve only a single type of polygon, namely honeycomb $=(6^3)$,
square $= (4^4)$, and triangular $=(3^6)$, abbreviated as (hc), (sq),
and (tri), respectively. The other Archimedean lattices are
heteropolygonal, i.e., they involve more than a single type of
polygon. Examples are $(4 \cdot 8 \cdot 8)$, $(3 \cdot 6 \cdot 3 \cdot
6)$ (often called ``kagom\'e'' (kag)), $(3^3 \cdot 4^2)$, and $(3
\cdot 3 \cdot 4 \cdot 3 \cdot 4)$.  With appropriate boundary
conditions, a finite section of an Archimedean lattice is a
$\Delta$-regular graph. The $\Delta$ values for the Archimedean
lattices range from 3 to 6.  Our upper bounds were denoted
$\phi_u(\Lambda)$ and $\sigma_u(\Lambda)$, where the subscript $u$
stands for ``upper''.  Our bounds are, to our knowledge, the best
upper bounds on $\phi(\Lambda)$ and $\sigma(\Lambda)$ for these
lattices. In \cite{sfc} we also calculated lower bounds on these
exponential growth constants, which are very close to the respective
upper bounds, with fractional differences ranging from $10^{-4}$ to
$10^{-2}$. This property enabled us to calculate quite accurate
approximate values of $\phi(\Lambda)$ and $\sigma(\Lambda)$ for these
lattices, which we denote here simply as $\phi(\Lambda)$ and
$\sigma(\Lambda)$.  For each lattice $\Lambda$, our method made use of 
calculations of lower and upper bounds on $\phi$ and $\sigma$
for a sequence of infinite-width lattice strips of increasing widths.
Our approximate values were determined conservatively as the average of
our lower and upper bounds on the widest strips with periodic transverse
boundary conditions (to minimize finite-width effects). 
As we noted in \cite{sfc}, our upper bounds and
approximate values for $\phi(\Lambda)$ and $\sigma(\Lambda)$ are
monotonically increasing functions of the vertex degree $\Delta$ for
the Archimedean lattices that we studied.  In the following, we focus on
our results for $\phi(\Lambda)$ in view of new general bounds in
\cite{bcl}.  For reference, we list our
upper bounds and values for $\phi(\Lambda)$ from \cite{sfc} in Table
\ref{phi_versus_upperbound_table}, together with the ratio of our
(central) value of $\phi(\Lambda)$ divided by our upper bound for
each $\Lambda$, namely
\beq
R_{\phi}(\Lambda) = \frac{\phi(\Lambda)}{\phi_u(\Lambda)} \ . 
\label{ratio_phi}
\eeq
The fact that the ratios $R_{\phi}(\Lambda)$ for these lattices
are very close to unity shows how close our upper bounds are to
being sharp. We found that our upper bounds on $\phi(\Lambda)$ and
$\sigma(\Lambda)$ approach limiting values more rapidly than our lower
bounds, so that the true values of $R_\phi(\Lambda)$ are expected to be
even closer to unity than the values listed in Table
\ref{phi_versus_upperbound_table}, i.e., the respective upper bounds are
even closer to being sharp for these lattices.

Let ${\cal G}_\Delta$ denote the set of all $\Delta$-regular
$n$-vertex graphs. Recently, in Ref. \cite{bcl}, Bob\'enyi,
Csikv\'ari, and Luo (BCL) presented upper bounds on the supremum over
all $G \in {\cal G}_\Delta$ of the quantity $[N_{SF}(G)]^{1/n(G)}$,
\beq
f_\Delta = {\rm sup}_{G \in {\cal G}_\Delta} [N_{SF}(G)]^{1/n(G)} \ .
\label{fd}
\eeq
The upper bounds reported in \cite{bcl} do not depend on $n(G)$, so
they also apply in the $n(G) \to \infty$ limit, yielding upper bounds
on $\phi$, which we denote as
\beq
\phi_{u,BCLi}(\Delta) = \lim_{n(G) \to \infty}
    {\rm sup}_{G \in {\cal G}_\Delta} [N_{SF}(G)]^{1/n(G)} \ ,
\label{phibcli}
\eeq
where the subscript $i$ will label the specific BCL bounds.  It is of
considerable interest to compare the BCL upper bounds with our upper
bounds and values for the Archimedean lattices that we considered in
\cite{sfc}.  We perform this comparison in the present paper.


\section{Comparison of Upper Bounds on $\phi(\Lambda)$}
\label{upper_bound_section} 

We first recall a general upper bound for any set of spanning
subgraphs, including spanning trees, spanning forests, and connected
spanning subgraphs.  In the construction of a spanning subgraph, there
is choice for each edge of $G$, namely whether it is present or
absent. Since this is a two-fold choice for each edge, it follows that
the number of spanning subgraphs of $G$ is 
\beq
N_{SSG}(G)=2^{e(G)} \ .
\label{nssg}
\eeq
This is an upper bound for any specific subclass of spanning subgraphs. 
Hence, in particular, for spanning forests, 
\beq
N_{SF}(G) \le N_{SSG}(G) \ .
\label{nsf_le_nssg}
\eeq
Since for $\Delta$-regular graphs $G \in {\cal G}_\Delta$, 
\beq
e(G) = \frac{n(G)\Delta}{2} \ ,
\label{erel}
\eeq
it follows that for these $\Delta$-regular graphs,
\beq
f_\Delta \le 2^{\Delta/2} \ ,
\label{fd_expbound}
\eeq
and hence, in the limit $n(G) \to \infty$,
\beq
\phi(\{ {\cal G}_\Delta \}) \le 2^{\Delta/2} \ .
\label{phi_expbound}
\eeq

Before comparing our upper bounds on $\phi(\Lambda)$ to the $n(G) \to
\infty$ limits of upper bounds recently derived in \cite{bcl}, we
mention some previous bounds.  After early work \cite{mw},
Ref. \cite{cmnn} obtained the upper limit
\beq
\phi(sq) \le 3.7410018 \ 
\label{cmnn_bound}
\eeq
Before our work in \cite{sfc}, the best upper bound on $\phi(sq)$ was from
Mani, in Ref. \cite{mani2012}, namely
\beq
\phi(sq) \le 3.705603 \ .
\eeq
In \cite{sfc} we derived the upper bound 
\beq
\phi(sq) \le 3.699659 \ .
\label{phiup_sfc}
\eeq
As we noted, our upper bounds on $\phi(\Lambda)$ for this and the
other Archimedean lattices that we studied are, to our knowledge,
the best upper bounds on $\phi(\Lambda)$ for these lattices. Our
results in \cite{sfc} were part of a general program of calculating
bound on, and values of, exponential growth constants for
various classes of subgraphs on Archimedean lattices \cite{ka3,aca,ac}.

A first upper bound proved in \cite{bcl} is
\beq
N_{SF}(G) \le \prod_{v_i \in V} (\Delta_{v_i}+1) \ .
\label{nsfbound1}
\eeq
The special case of this bound for a $\Delta$-regular graph $G$ is  
\beq
f_\Delta \le \Delta + 1 \ , 
\label{fd_linbound}
\eeq
which also applies in the limit as $n(G) \to \infty$ as
\beq
\phi(\{ {\cal G}_\Delta \}) \le \Delta+1 \ .
\label{phi_linbound}
\eeq
We denote the right-hand side of (\ref{phi_linbound}) as  
$\phi_{u,\Delta,BCL1}(\Delta) = \Delta + 1$. 
Generalizing $\Delta$ from positive integral values to positive real values,
we find that the upper bound (\ref{phi_expbound}) is more stringent than
(\ref{phi_linbound}) if $\Delta < 5.3197$.

A second upper bound for $\Delta$-regular graphs discussed in \cite{bcl} is
\beq
f_\Delta \le f_{u,BCL2}(\Delta) \ , 
\label{fd_bcl2}
\eeq
where 
\beq
f_{u,BCL2}(\Delta) = \bigg ( \frac{\Delta+1}{\eta(\Delta)} \bigg )
\bigg ( \frac{\Delta-1}{\Delta -\eta(\Delta)} \bigg )^{\frac{\Delta-2}{2}} ,
\label{fd_bcl2b}
\eeq
with 
\beq
\eta(\Delta) = \frac{(\Delta+1)(\Delta+1-\sqrt{\Delta^2 -2\Delta+5} \ )}
    {2(\Delta -1)} \ .
\label{eta}
\eeq
(See also \cite{ks} for related work.) 
Since this bound applies uniformly for any $n(G)$, it also applies to the
limit as $n(G) \to \infty$:
\beq
\phi(\Lambda) \le \phi_{u,BCL2}(\Delta) \ , 
\label{phiup_bcl2}
\eeq
where $\phi_{u,BCL2}(\Delta) = f_{u,BCL2}(\Delta)$. 
We list below the analytic expressions of the upper bound
$\phi_{u,BCL2}(\Delta)$ and the corresponding numerical values (given to
the indicated number of significant figures) for the values of 
$\Delta$ that are relevant for comparison with our bounds: 
\beq
\phi_{u,BCL2}(2) = \frac{3+\sqrt{5}}{2} = 2.618034
\label{phi_bcl2_delta2}
\eeq
\beq
\phi_{u,BCL2}(3) =2 \bigg ( \frac{11+8\sqrt{2}}{7} \bigg )^{1/2} = 3.5708109 
\label{phi_bcl2_delta3}
\eeq
\beq
\phi_{u,BCL2}(4) = \frac{35+13\sqrt{13}}{18} = 4.5484537 
\label{phi_bcl2_delta4}
\eeq
\beq
\phi_{u,BCL2}(5)=\frac{4\sqrt{2} \, (27+7\sqrt{5})\sqrt{1+3\sqrt{5}}}{121}
  =5.5361833
\label{phi_bcl2_delta5}
\eeq
and
\beq
\phi_{u,BCL2}(6) = \frac{791+58\sqrt{29}}{169} = 6.52863644 \ .
\label{phi_bcl2_delta6}
\eeq

An upper bound given in Ref. \cite{bcl} for 4-regular graphs ${\cal G}_4$,
which we label BCL3, is:
\beq
f_4 \le 3.994 \ , 
\label{f4}
\eeq
and again, since this is independent of $n(G)$ for $G \in {\cal G}_4$, it
implies, in the $n(G) \to \infty$ limit, the upper bound
$\phi(\{ {\cal G}_4 \}) \le \phi_{u,BCL3}(4)$, where 
\beq
\phi_{u,BCL3}(4) = 3.994  \ . 
\label{bcl_theorem1p5}
\eeq

Finally, Ref. \cite{bcl} presented slightly stronger upper bounds
on $f_\Delta$ for $\Delta$-regular graphs ${\cal G}_\Delta$ with
$\Delta$ in the interval $5 \le \Delta \le 9$. As before, their upper
bound for each $\Delta$ makes no reference to $n(G)$, so that it
implies the same bound for the $n(G) \to \infty$ limit of a
$\Delta$-regular graph for $\Delta$ in this interval $5 \le \Delta \le
9$:
\beq
\phi(\Lambda) \le \phi_{u,BCL4}(\Delta) \quad {\rm for} \ 5 \le \Delta \le 9
\ .
\label{phiup_bcl4}
\eeq
In the two cases with $\Delta=5$ and $\Delta=6$ relevant for
Archimedean lattices, these upper bounds are $\phi_{u,BCL4}(5)=5.1965$
and $\phi_{u,BCL4}(6) = 6.3367$ \cite{bcl}.  We list the numerical
values of the BCL upper bounds on $\phi(\Lambda)$ in Table
\ref{phiup_cs_versus_phiup_ks_bcl_table} for the $\Delta$ values that
occur for Archimedean lattices, namely $\Delta=3,4,5,6$.  As is
evident, for a given $\Delta$, these BCL upper bounds are all less
stringent than the upper bounds that we derived on $\phi(\Lambda)$ for
the $\Delta$-regular Archimedean lattices $\Lambda$ in \cite{sfc}. This
comparison is made for the full range of $\Delta$ values on
Archimedean lattices, namely $\Delta=3, \ 4, \ 5, \ 6$. For the 
case of $\Delta=2$, we recall the elementary result that 
$N_{SF}(C_n)=2^n-1$ and hence, in the $n \to \infty$ limit,
$\phi(C_\infty)=2$. This saturates the upper bound
(\ref{phi_expbound}), but is less than the upper bounds
$\phi_{u,BCL1}(2)=3$ and $\phi_{u,BCL2}(2)=(3+\sqrt{5})/2=2.6180$.

Although the BCL bounds on $\phi$ depend only on $\Delta$, our upper
bounds $\phi_u(\Lambda)$ and also the values of $\phi(\Lambda)$ that
we calculated in \cite{sfc} enable us to investigate the dependence on
girth $g(\Lambda)$ for lattices with the same value of $\Delta$.  As
one can see from Table \ref{phi_versus_upperbound_table}, two such
comparisons can be made for Archimedean lattices: (i) the $(6^3)=$hc
and $(4 \cdot 8^2)$ lattices both have $\Delta=3$, but the girth
$g(hc)=6$ is larger than the girth $g((4 \cdot 8^2))=4$, and
$\phi(hc)$ is slightly larger than $\phi((4 \cdot 8^2))$; (ii) the
square and kagom\'e lattices both have $\Delta=4$, but the girth
$g(sq)=4$ is larger than the girth $g(kag)=3$, and $\phi(sq)$ is
slightly larger than $\phi(kag)$.  Note that the uncertainties in our
numerical determination of the values of $\phi(\Lambda)$ for the
Archimedean lattices are sufficiently small that they are negligible
for these comparisons. 

In earlier work preceding \cite{ac,aca,sfc}, we had calculated 
values of exponential growth constants for a variety of families of
lattice strip graphs of various lattices with a range of finite widths
and with arbitrarily great length (e.g., \cite{a}-\cite{sdg}).  We
found that for a given type of lattice strip graph, in the
infinite-length limit, $\phi$ is a monotonically increasing function
of the strip width.  For the smallest widths, one obtains simple
algebraic expressions for these exponential growth constants. For
example, $\phi$ has the following values for the infinite-length
limits of the given strips: 
$\phi=1+\sqrt{3}=2.732$ and $\phi=(\sqrt{7}+\sqrt{15})/2=3.259$ for
the transverse width $L_t=2$ strips of the square lattice with free
(F) and periodic (P) tranverse boundary conditions (BC$_t$);
$\phi=2+\sqrt{2}=3.414$ and $\phi=[(23+\sqrt{505})/2]^{1/2}=4.768$
for $L_t=2$ strips of the
triangular lattice with F and P BC$_t$, etc.

While our upper bounds $\phi_u(\Lambda)$ are better than the BCL upper
bounds in \cite{bcl} for the Archimedean lattices $\Lambda$ that we
studied, the BCL upper bounds are still valuable, since they apply for any
$\Delta$-regular graphs, not just Archimedean lattices.  
In future work it would be of interest to search for $\Delta$-regular
families of graphs for which
$f_\Delta$ and/or $\lim_{n(G) \to \infty} f_\Delta$ lie closer to the
upper bounds given in \cite{bcl}.  


\begin{acknowledgments}

This research was supported in part by the Taiwan Ministry of Science and
Technology grant MOST 109-2112-M-006-008 (S.-C.C.) and by
the U.S. National Science Foundation grant No. NSF-PHY-1915093 (R.S.).

\end{acknowledgments}


\begin{appendix}

\section{Some Background from Graph Theory}
\label{graphtheory}

In this appendix we briefly review some background from graph theory,
in particular, a connection of $N_{SF}(G)$ and
$N_{CSSG}(G)$ with evaluations of the Tutte polynomial.  As in the text, let
$G=(V,E)$ be a graph defined by its vertex and edge sets $V$ and
$E$. Further, let $n=n(G)=|V|$, $e(G)=|E|$, $k(G)$, and $c(G)$ denote
the numbers of vertices, edges, connected components, and linearly
independent cycles in $G$, respectively. The Tutte polynomial of a
graph $G$, denoted $T(G,x,y)$, is defined as
\beq
T(G,x,y) = \sum_{G' \subseteq G} (x-1)^{k(G')-k(G)} \, (y-1)^{c(G')} \ ,
\label{t}
\eeq
where $G'$ is a spanning subgraph of $G$ (see, e.g.,
\cite{tutte67,graphtheory}).  The numbers of spanning forests and
connected spanning subgraphs are evaluations of the Tutte polynomial:
\beq
N_{SF}(G) = T(G,2,1)
\label{nsf_eq_tx2y1}
\eeq
and
\beq
N_{CSSG}(G) = T(G,1,2) \ .
\label{nccsg_eq_tx1y2}
\eeq
For a general graph $G$, the calculation
of $N_{SF}(G)$ and $N_{CSSG}(G)$  are $\sharp$ P hard \cite{jvw}.
This is why it is useful to have bounds on these quantities, and also on
the corresponding exponential growth constants.

A remark is in order here concerning graphs with loops. Recall that 
a loop is an edge that connects a vertex back to itself.  The reason
that we restrict to loopless graphs in our work is that if one allows
loops, then one loses a connection between the vertex degree of a
$\Delta$-regular graph $G$ and $\phi(\{G\})$.  This can be illustrated
in the simple case of the circuit graph $C_n$, which is
$\Delta$-regular with $\Delta=2$.  One has $T(C_n,x,y) =
y+\sum_{j=1}^{n-1} x^j$, so that, in the limit $n \to \infty$,
$\phi(\{ C \})=2$. Now let us attach $m$ loops ($\ell$) to each vertex of
$C_n$.  We denote the resultant graph as $C_{n,m\ell}$.  This is again a
$\Delta$-regular graph with vertex degree $\Delta=2(1+m)$. The Tutte polynomial
is
\beq
T(C_{n,m\ell},x,y)=y^{mn}T(C_n,x,y)=y^{mn}\bigg (
y+\sum_{j=1}^{n-1} x^j \bigg ) \ .
\label{tcc}
\eeq
Hence,
\beq
N_{SF}(C_{n,m\ell})=T(C_{n,m\ell},2,1)=T(C_n,2,1)=N_{SF}(C_n)
\label{nsfcc}
\eeq
and, in the limit $n \to \infty$, the corresponding values of $\phi$
are the same for the $C_n$ and $C_{n,m\ell}$ families of graphs,
although the vertex degrees are different for these families.  Thus,
if one were to allow modifications of Archimedean lattices with loops,
one would lose the informative connection between the vertex degree and
the value of $\phi(\Lambda)$.

\end{appendix}



\newpage


\begin{table*}
  \caption{\footnotesize{For each Archimedean lattice $\Lambda$, this table
      lists the value of $\phi(\Lambda)$ and the upper
      bound, $\phi_u(\Lambda)$, both from Ref. \cite{sfc},
      together with the ratio $R_{\phi}(\Lambda) =
      \phi(\Lambda)/\phi_u(\Lambda)$.  The lattices are listed in
      order of increasing vertex degree $\Delta(\Lambda)$.}}
  \begin{center}
  \begin{ruledtabular}
\begin{tabular}{|c|c|c|c|c|c|} 
$\Lambda$ & $\Delta(\Lambda)$ & $g(\Lambda)$ & $\phi(\Lambda)$ & 
$\phi_u(\Lambda)$ & $R_{\phi}(\Lambda)$ \\
  \hline
$(4 \cdot 8^2)$ & 3 & 4 & $2.77931 \pm 0.00018$ & 2.779486 & 0.99994 \\ \hline 
$(6^3)=$ hc     & 3 & 6 & $2.80428 \pm 0.00050$ & 2.804781 & 0.99982 \\ \hline 
$(3\cdot 6\cdot 3\cdot 6)$
                & 4 & 3 & $3.602 \pm 0.012$     & 3.614045 & 0.99667 \\ \hline 
$(4^4)=$ sq     & 4 & 4 & $3.687 \pm 0.012$     & 3.699659 & 0.99658 \\ \hline 
$(3^3 \cdot 4^2)$& 5 &3 & $4.530 \pm 0.024$     & 4.553665 & 0.99480 \\ \hline 
$(3^2 \cdot 4 \cdot 3 \cdot 4)$
                 & 5 &3 & $4.503 \pm 0.065$    & 4.568231 & 0.98572 \\ \hline
$(3^6)=$ tri     & 6 &3 & $5.444 \pm 0.051$    & 5.494840 & 0.99075 \\ 
\end{tabular}
\end{ruledtabular}
\end{center}
\label{phi_versus_upperbound_table}
\end{table*}


\begin{table*}
  \caption{\footnotesize{Comparison of upper bounds on $\phi(\Lambda)$
      for Archimedean lattices $\Lambda$. The most stringent upper
      bounds on $\phi(\Lambda)$ are those from Ref. \cite{sfc},
      denoted $\phi_u(\Lambda)$. The table also lists the general
      upper bound $2^{\Delta/2}$ and, where applicable, the upper
      bounds $\phi_{u,BCLi}(\Delta)$, $i=1,2,3,4$ from \cite{bcl}.
      The BCL3 bound applies for $\Delta=4$, while the BCL4 bound
      applies for $\Delta=5, 6$. For lattices where a given BCLi bound
      is not applicable, we denote this by a dash.}}
\begin{center}
\begin{ruledtabular}
\begin{tabular}{|c|c|c|c|c|c|c|c|} 
  $\Lambda$ & $\Delta(\Lambda)$ & $g(\Lambda)$ &  $\phi_u(\Lambda)$ &
  $2^{\Delta/2}$ & $\phi_{u,BCL1}(\Delta)$ & $\phi_{u,BCL2}(\Delta)$ &
  $\phi_{u,BCL3,4}(\Delta)$ \\
  \hline
  $(4 \cdot 8^2)$&3 & 4     & 2.779486 & 2.82843 & 4 & 3.57081 & $-$ \\ \hline 
  $(6^3)=$ hc   & 3 & 6     & 2.804781 & 2.82843 & 4 & 3.57081 & $-$ \\ \hline 
  $(3\cdot 6\cdot 3\cdot 6)$
                & 4 & 3     & 3.614045 & 4       & 5 & 4.54845 & 3.994 \\ \hline 
  $(4^4)=$ sq   & 4 & 4     & 3.699659 & 4       & 5 & 4.54845 & 3.994 \\ \hline 
  $(3^3 \cdot 4^2)$&5 & 3   & 4.553665 & 5.65685 & 6 & 5.53618 & 5.1965 \\ \hline 
  $(3^2 \cdot 4 \cdot 3 \cdot 4)$
                & 5 & 3     & 4.568231 & 5.65685 & 6 & 5.53618 & 5.1965 \\ \hline
  $(3^6)=$ tri  & 6 & 3     & 5.494840 & 8       & 7 & 6.52864 & 6.3367  \\ 
\end{tabular}
\end{ruledtabular}
\end{center}
\label{phiup_cs_versus_phiup_ks_bcl_table}
\end{table*}


\end{document}